\documentclass[12pt,twoside]{article}
\usepackage{amsmath, amssymb, amsthm,graphicx}
\font\fourteenb=cmb10 at 14pt
\begin{document}
\vspace*{-1.0cm}\noindent \copyright
 Journal of Technical University at Plovdiv\\[-0.0mm]\
\ Fundamental Sciences and Applications, Vol. 1, 1995\\[-0.0mm]
\textit{Series A-Pure and Applied Mathematics}\\[-0.0mm]
\ Bulgaria, ISSN 1310-8271\\[+1.2cm]
\font\fourteenb=cmb10 at 12pt
\begin{center}

{\bf \LARGE Properties of special classes of analytic functions
which are multipliers of the Cauchy-Stieltjes integrals
   \\ \ \\ \large Peyo  Stoilov, Milena Racheva }
\end{center}


\begin{abstract}
Special classes of analytic functions, denoting by $K$  and  $J$
arc considered in this paper. It is proved, that $f\in J$   if and
only if   $f'\in K$  ;  $f\in J$  yildse $f\circ \phi _{a} \in J$
, where  $\phi _{a} (z)=z+a/1+za,{\kern 1pt} {\kern 1pt} {\kern
1pt} {\kern 1pt} {\kern 1pt} {\kern 1pt} {\kern 1pt} {\kern 1pt}
-1<a<1;$  $J\subset F_{0} \subset BMOA$ . The classes  $J$  and
$K$   were introduced in [4].
\end{abstract}

     \section{Introduction}

Let $D$ denote the unit disk in the complex plane and $T$ the unit circle.
Let $M$ be the Banach space of all complex-valued Borel measures on $T$ with
the usual variation norm.

For $\alpha \ge 0$ , let $F_{\alpha } $ denote the family of the analytic
functions $f$ , for which there exists $\mu \in M$ such that
\begin{equation*}
\displaystyle f(z)=\int _{T}\frac{1}{(1-\overline{\xi }z)^{\alpha } } {\kern %
1pt} {\kern 1pt} {\kern 1pt} {\kern 1pt} {\kern 1pt} d\mu (\xi ),{\kern 1pt}
{\kern 1pt} {\kern 1pt} {\kern 1pt} {\kern 1pt} {\kern 1pt} {\kern 1pt} {%
\kern 1pt} {\kern 1pt} {\kern 1pt} {\kern 1pt} {\kern 1pt} {\kern 1pt} {%
\kern 1pt} {\kern 1pt} {\kern 1pt} {\kern 1pt} {\kern 1pt} {\kern 1pt} {%
\kern 1pt} {\kern 1pt} {\kern 1pt} {\kern 1pt} {\kern 1pt} {\kern 1pt} {%
\kern 1pt} {\kern 1pt} {\kern 1pt} {\kern 1pt} {\kern 1pt} {\kern 1pt} {%
\kern 1pt} {\kern 1pt} {\kern 1pt} {\kern 1pt} {\kern 1pt} {\kern 1pt} {%
\kern 1pt} {\kern 1pt} {\kern 1pt} {\kern 1pt} {\kern 1pt} {\kern 1pt} {%
\kern 1pt} {\kern 1pt} {\kern 1pt} {\kern 1pt} {\kern 1pt} {\kern 1pt} {%
\kern 1pt} {\kern 1pt} {\kern 1pt} {\kern 1pt} {\kern 1pt} {\kern 1pt} {%
\kern 1pt} {\kern 1pt} {\kern 1pt} {\kern 1pt} {\kern 1pt} {\kern 1pt}
\alpha >0,
\end{equation*}
\begin{equation*}
\displaystyle f(z)=f(0)+\int _{T}\log \left(\frac{1}{1-\overline{\xi }z}
\right) {\kern 1pt} {\kern 1pt} {\kern 1pt} {\kern 1pt} {\kern 1pt} d\mu
(\xi ),{\kern 1pt} {\kern
1pt} {\kern 1pt} {\kern 1pt} {\kern 1pt} {\kern 1pt} {\kern 1pt} {\kern 1pt}
{\kern 1pt} {\kern 1pt} {\kern 1pt} {\kern 1pt} {\kern
1pt} {\kern 1pt} {\kern 1pt} {\kern 1pt} {\kern 1pt} {\kern 1pt} {\kern 1pt}
{\kern 1pt} {\kern 1pt} {\kern 1pt} {\kern 1pt} {\kern
1pt} {\kern 1pt} {\kern 1pt} \alpha =0.
\end{equation*}

Let $\mu \in M$ . In [4] were introduced the functions $J_{\mu } $ and $%
K_{\mu } $ by

\begin{equation*}
J_{\mu } (z)=\int _{T}\log \left(\frac{1-\xi z}{1-\overline{\xi }z} \right) {%
\kern 1pt} {\kern 1pt} {\kern 1pt} {\kern 1pt} {\kern 1pt} d\mu (\xi ),{%
\kern 1pt} {\kern 1pt} {\kern 1pt} {\kern 1pt} {\kern 1pt} {\kern 1pt} {%
\kern 1pt} {\kern 1pt} {\kern 1pt} {\kern 1pt} {\kern 1pt} {\kern 1pt} {%
\kern 1pt} {\kern 1pt} {\kern 1pt} {\kern 1pt} {\kern 1pt} {\kern 1pt} {%
\kern 1pt} {\kern 1pt} {\kern 1pt} {\kern 1pt} {\kern 1pt} {\kern 1pt} {%
\kern 1pt}
\end{equation*}

\begin{equation*}
K_{\mu } (z)=\int _{T}\left(\frac{1}{\xi -z} -\frac{\xi }{1-\xi z} \right) {%
\kern 1pt} {\kern 1pt} {\kern 1pt} {\kern 1pt} {\kern 1pt} d\mu (\xi )=\int
_{T}\frac{1-\xi ^{2} }{(\xi -z)(1-\xi z)} {\kern 1pt} {\kern 1pt} {\kern 1pt}
{\kern 1pt} {\kern 1pt} d\mu (\xi ).
\end{equation*}

In this paper some properties of the classes $J$ and $K$ are given, where

\begin{equation*}
J=\left\{f:{\kern 1pt} {\kern 1pt} {\kern 1pt} {\kern 1pt} {\kern 1pt} {%
\kern 1pt} {\kern 1pt} {\kern 1pt} {\kern 1pt} f=J_{\mu } ,{\kern 1pt} {%
\kern 1pt} {\kern 1pt} {\kern 1pt} {\kern 1pt} {\kern 1pt} {\kern 1pt} {%
\kern 1pt} {\kern 1pt} {\kern 1pt} {\kern 1pt} \mu \in M\right\},
\end{equation*}

\begin{equation*}
K=\left\{f:{\kern 1pt} {\kern 1pt} {\kern 1pt} {\kern 1pt} {\kern 1pt} {%
\kern 1pt} {\kern 1pt} {\kern 1pt} {\kern 1pt} f=K_{\mu } ,{\kern 1pt} {%
\kern 1pt} {\kern 1pt} {\kern 1pt} {\kern 1pt} {\kern 1pt} {\kern 1pt} {%
\kern 1pt} {\kern 1pt} {\kern 1pt} {\kern 1pt} \mu \in M\right\}.
\end{equation*}

\section{Main results}

Let
\begin{equation*}
\phi _{a} (z)=\frac{z+a}{1+za} ,{\kern 1pt} {\kern 1pt} {\kern 1pt} {\kern %
1pt} {\kern 1pt} {\kern 1pt} {\kern 1pt} {\kern 1pt}  -1<a<1.
\end{equation*}

\textbf{Lemma 1.} \textit{If } $f\in K$ \textit{then } $(f\circ \phi _{a} ){%
\kern 1pt} {\kern 1pt} {\kern 1pt} {\kern 1pt} \phi ^{\prime }_{a} {\kern 1pt%
} {\kern 1pt} {\kern 1pt} {\kern 1pt} {\kern 1pt} {\kern 1pt} {\kern 1pt} {%
\kern 1pt} {\kern 1pt} \in {\kern 1pt} {\kern 1pt} {\kern 1pt} {\kern 1pt} {%
\kern 1pt} {\kern 1pt} K.$

\textit{Proof.} Suppose that $f=K_{\mu } .$ Then

\begin{equation*}
(f\circ \phi _{a} ){\kern 1pt} =\int _{T}\left(\frac{1}{\xi -\frac{z+a}{1+za}
} -\frac{1}{1-\xi \frac{z+a}{1+za} } \right) {\kern 1pt} {\kern 1pt} {\kern %
1pt} {\kern 1pt} {\kern 1pt} d\mu (\xi )=
\end{equation*}

\begin{equation*}
=(1+za)\int _{T}\left(\frac{1}{\xi (1+za)-z+a} -\frac{1}{(1+za)-\xi (z+a)}
\right) {\kern 1pt} {\kern 1pt} {\kern 1pt} {\kern 1pt} {\kern 1pt} d\mu
(\xi )=
\end{equation*}

\begin{equation*}
=(1+za)\int _{T}{\kern 1pt} {\kern 1pt} {\kern 1pt} \frac{1}{1-\xi a} {\kern %
1pt} {\kern 1pt} {\kern 1pt} {\kern 1pt} {\kern 1pt} {\kern 1pt} \left(\frac{%
1}{\frac{\xi -a}{1-\xi a} -z} -\frac{\xi }{1-\frac{\xi -a}{1-\xi a} {\kern %
1pt} {\kern 1pt} {\kern 1pt} {\kern 1pt} z} \right) {\kern 1pt} {\kern 1pt} {%
\kern 1pt} {\kern 1pt} {\kern 1pt} d\mu (\xi ).
\end{equation*}
Since $\left|a\right|<1$ , the function
\begin{equation*}
\phi _{a} ^{-1} (z)=\frac{z-a}{1-za}
\end{equation*}
maps $T$ one-to-one on $T$ .  For each Borel set $E\subset T,$ let $\nu
(E)=\mu (\phi _{a} (E)).$ Then $\nu \in M$ and we can write
\begin{equation*}
f\circ \phi _{a} {\kern 1pt} =(1+za)\int _{T}{\kern 1pt} {\kern 1pt} {\kern %
1pt} \frac{1}{1-\phi _{a} (\xi ){\kern 1pt} {\kern 1pt} {\kern 1pt} a} {%
\kern 1pt} {\kern 1pt} {\kern 1pt} {\kern 1pt} {\kern 1pt} {\kern 1pt} \left(%
\frac{1}{\xi -z} -\frac{\phi _{a} (\xi ){\kern 1pt} }{1-\xi {\kern 1pt} {%
\kern 1pt} {\kern 1pt} z} \right) {\kern 1pt} {\kern 1pt} {\kern 1pt} {\kern %
1pt} {\kern 1pt} d\nu (\xi )=
\end{equation*}
\begin{equation*}
=\frac{(1+za)^{2} }{1-a^{2} } \int _{T}\frac{1-\xi ^{2} }{(\xi -z)(1-\xi z)}
{\kern 1pt} {\kern 1pt} {\kern 1pt} {\kern 1pt} {\kern 1pt} d\nu (\xi ){%
\kern 1pt} {\kern 1pt} {\kern 1pt} {\kern 1pt} {\kern 1pt} {\kern 1pt} ={%
\kern 1pt} {\kern 1pt} {\kern 1pt} {\kern 1pt} {\kern 1pt} {\kern 1pt} \frac{%
1}{\phi ^{\prime }_{a} (z)} {\kern 1pt} {\kern 1pt} {\kern 1pt} {\kern 1pt}
\cdot q(z),
\end{equation*}
where $q=K_{\nu } {\kern 1pt} {\kern 1pt} {\kern 1pt} \in K.$  Consequently
\begin{equation*}
(f\circ \phi _{a} ){\kern 1pt} {\kern 1pt} {\kern 1pt} {\kern 1pt} \phi
^{\prime }_{a} {\kern 1pt} {\kern 1pt} {\kern 1pt} {\kern 1pt} =q(z){\kern %
1pt} {\kern 1pt} {\kern 1pt} {\kern 1pt} {\kern 1pt} \in {\kern 1pt} {\kern %
1pt} {\kern 1pt} {\kern 1pt} {\kern 1pt} {\kern 1pt} K.
\end{equation*}

\textbf{Theorem 1.}

\textit{a) } $f(z)-f(0)=J_{\mu } (z){\kern 1pt} {\kern 1pt} {\kern 1pt} {%
\kern 1pt} {\kern 1pt} $ \textit{if }\textit{and }\textit{only }\textit{if }
$f^{\prime }(z)=K_{\mu } (z){\kern 1pt} {\kern 1pt} ; $

\textit{b}\textit{)} $J\subset F_{\alpha } $ \textit{for every} $\alpha \ge
0;$

\textit{c) }\textit{If } $f\in J$ \textit{then } $f\circ \phi _{a} \in J.$

\

\textit{Proof.} Suppose that $f(z)-f(0)=J_{\mu } (z){\kern 1pt} {\kern 1pt} {%
\kern 1pt} {\kern 1pt} {\kern 1pt} .$

\

Then $f$ has the representation
\begin{equation*}
f(z)-f(0)=\int _{T}\log \left(\frac{1-\xi z}{1-\overline{\xi }z} \right) {%
\kern 1pt} {\kern 1pt} {\kern 1pt} {\kern 1pt} {\kern 1pt} d\mu (\xi ),{%
\kern 1pt} {\kern 1pt} {\kern 1pt} {\kern 1pt}
\end{equation*}
with $\mu \in M.$  Then
\begin{equation*}
f^{\prime }(z)=\int _{T}\left(\frac{-\xi }{1-\xi z} +\frac{\xi }{1-\overline{%
\xi }z} \right) {\kern 1pt} {\kern 1pt} {\kern 1pt} {\kern 1pt} {\kern 1pt}
d\mu (\xi )=\int _{T}\left(\frac{-\xi }{1-\xi z} +\frac{1}{\xi -z} \right) {%
\kern 1pt} {\kern 1pt} {\kern 1pt} {\kern 1pt} {\kern 1pt} d\mu (\xi )=K\mu .
\end{equation*}
Conversety, suppose that $f^{\prime }(z)=K_{\mu } (z){\kern 1pt} {\kern 1pt}
,{\kern 1pt} {\kern 1pt} {\kern 1pt} $ i.e

\begin{equation*}
f^{\prime }(z)=\int _{T}\left(\frac{1}{\xi -z} -\frac{\xi }{1-\xi z} \right)
{\kern 1pt} {\kern 1pt} {\kern 1pt} {\kern 1pt} {\kern 1pt} d\mu (\xi ), {%
\kern 1pt} {\kern 1pt}{\kern 1pt} {\kern 1pt}\mu \in M.
\end{equation*}
Then
\begin{equation*}
f(z)=\int _{0}^{z}f^{\prime }(t) {\kern 1pt} {\kern 1pt} {\kern 1pt}
dt+f(0)=.
\end{equation*}

\begin{equation*}
=\int _{0}^{z}\left(\int _{T}\left(\frac{1}{\xi -t} -\frac{\xi }{1-\xi t}
\right) {\kern 1pt} {\kern 1pt} {\kern 1pt} {\kern 1pt} {\kern 1pt} d\mu
(\xi )\right) {\kern 1pt} {\kern 1pt} dt+f(0)=
\end{equation*}

\begin{equation*}
=\int _{T}{\kern 1pt} \left(\int _{0}^{z}\left(\frac{1}{\xi -t} -\frac{\xi }{%
1-\xi t} \right) {\kern 1pt} {\kern 1pt} {\kern 1pt} d{\kern 1pt} {\kern 1pt}
t\right){\kern 1pt} {\kern 1pt} {\kern 1pt} {\kern 1pt} {\kern 1pt} {\kern %
1pt} d\mu (\xi )+f(0){\kern 1pt} {\kern 1pt} {\kern 1pt} {\kern 1pt} {\kern %
1pt} {\kern 1pt} ={\kern 1pt} {\kern 1pt}
\end{equation*}

\begin{equation*}
=\int _{T}{\kern 1pt} \left(\log {\kern 1pt} {\kern 1pt} {\kern 1pt} \frac{%
1-\xi z}{\xi -z} -\log \frac{1}{\xi } \right){\kern 1pt} {\kern 1pt} {\kern %
1pt} {\kern 1pt} {\kern 1pt} {\kern 1pt} d\mu (\xi )+f(0){\kern 1pt} {\kern %
1pt} {\kern 1pt} {\kern 1pt} {\kern 1pt} {\kern 1pt} ={\kern 1pt} {\kern 1pt}
\end{equation*}

\begin{equation*}
=\int _{T}\log \left(\frac{1-\xi z}{1-\overline{\xi }z} \right) {\kern 1pt} {%
\kern 1pt} {\kern 1pt} {\kern 1pt} {\kern 1pt} d\mu (\xi ){\kern 1pt} {\kern %
1pt} {\kern 1pt} {\kern 1pt} {\kern 1pt} {\kern 1pt} +f(0){\kern 1pt} {\kern %
1pt} {\kern 1pt} {\kern 1pt} {\kern 1pt} =J_{\mu } (z)+f(0){\kern 1pt} {%
\kern 1pt} {\kern 1pt} .{\kern 1pt} {\kern 1pt} {\kern 1pt}
\end{equation*}

\

To prove \textit{b)}, assume that $f(z)=J_{\mu } (z){\kern 1pt} .{\kern 1pt}
{\kern 1pt} {\kern 1pt} {\kern 1pt} $

Then \textit{a)} implies $f^{\prime }(z)=K_{\mu } (z){\kern 1pt} {\kern 1pt}
.{\kern 1pt} {\kern 1pt} {\kern 1pt} $

Since

\begin{equation*}
f^{\prime }(z)=\int _{T}\left(\frac{1}{\xi -z} -\frac{\xi }{1-\xi z} \right)
{\kern 1pt} {\kern 1pt} {\kern 1pt} {\kern 1pt} {\kern 1pt} d\mu (\xi
)=-\int _{T}\frac{\xi }{1-\xi z} {\kern 1pt} {\kern 1pt} {\kern 1pt} {\kern %
1pt} {\kern 1pt} d\mu (\xi )-\int _{T}\frac{\xi }{1-\xi z} {\kern 1pt} {%
\kern 1pt} {\kern 1pt} {\kern 1pt} {\kern 1pt} d\mu (\overline{\xi }),
\end{equation*}

then $f^{\prime }(z)\in F_{1} (z){\kern 1pt} {\kern 1pt} .{\kern 1pt} {\kern %
1pt} {\kern 1pt} $

From $f^{\prime }(z)\in F_{1} (z){\kern 1pt} {\kern 1pt} {\kern 1pt} {\kern %
1pt} {\kern 1pt} $ it follows that $f(z)=J_{\mu } (z){\kern 1pt} {\kern 1pt}
{\kern 1pt} {\kern 1pt} \in {\kern 1pt} {\kern 1pt} {\kern 1pt} F_{0} {\kern %
1pt} {\kern 1pt} $ [3].

Theorem 3[3] implies that $F_{0} \subset F_{\alpha } $ for $\alpha >0$ ,

therefore $J_{\mu }(z){\kern1pt}\in {\kern1pt}{\kern1pt}{\kern1pt}F_{\alpha }%
{\kern1pt}$ for every $\alpha \geq 0.$

\

We shall prove \textit{c)}. If $f\in J$ then $f^{\prime }(z)\in K{\kern 1pt}
{\kern 1pt} $ and by Lemma 1
\begin{equation*}
(f^{\prime }\circ \phi _{a} ){\kern 1pt} {\kern 1pt} {\kern 1pt} {\kern 1pt}
\phi ^{\prime }_{a} {\kern 1pt} {\kern 1pt} {\kern 1pt} {\kern 1pt} {\kern %
1pt} {\kern 1pt} {\kern 1pt} {\kern 1pt} {\kern 1pt} \in {\kern 1pt} {\kern %
1pt} {\kern 1pt} {\kern 1pt} {\kern 1pt} {\kern 1pt} K.
\end{equation*}

Consequently $(f\circ \phi _{a} )^{\prime }\in K$ which yields $f\circ \phi
_{a} \in J.$ \

\

\textbf{Theorem 2.} If $f(z){\kern 1pt} {\kern 1pt} \in {\kern 1pt} {\kern %
1pt} {\kern 1pt} F_{0} {\kern 1pt} {\kern 1pt} $ then $f(z){\kern 1pt} {%
\kern 1pt} \in {\kern 1pt} {\kern 1pt} {\kern 1pt} BMOA.{\kern 1pt} {\kern %
1pt} $

\

\  \textit{Proof.} Let $f(z){\kern 1pt} {\kern 1pt} \in {\kern 1pt} {\kern %
1pt} {\kern 1pt} F_{0} {\kern 1pt} {\kern 1pt} .$ Then $f(z){\kern 1pt} {%
\kern 1pt} $ has the representation

\begin{equation*}
f(z)=\int _{T}\log \left(\frac{1}{1-\overline{\xi }z} \right) {\kern 1pt} {%
\kern 1pt} {\kern 1pt} {\kern 1pt} {\kern 1pt} d\mu (\xi )=\sum \nolimits
_{n=1}^{\infty }\frac{z^{n} }{n} \int _{T}\xi ^{n} {\kern 1pt} {\kern 1pt} {%
\kern 1pt} {\kern 1pt} d\mu (\xi )
\end{equation*}

and
\begin{equation*}
\left|\hat{f}(n){\kern 1pt} \right|{\kern 1pt} {\kern 1pt} {\kern 1pt}
=\left|\frac{1}{n} \int _{T}\xi ^{n} {\kern 1pt} {\kern 1pt} {\kern 1pt} {%
\kern 1pt} d\mu (\xi )\right|\le \frac{1}{n} \left\| \mu \right\| .
\end{equation*}

Now we shall make use of the fact that ${\kern 1pt} {\kern 1pt} {\kern 1pt} {%
\kern 1pt} BMOA\approx H^{1*} {\kern 1pt} $ [1]\textit{, as}

\begin{equation*}
<f,h>=\mathop{\lim }\limits_{r\to 1} \int _{T}\overline{f(r\xi )} {\kern 1pt}
{\kern 1pt} {\kern 1pt} {\kern 1pt} h(\xi ){\kern 1pt} {\kern 1pt} {\kern 1pt%
} d{\kern 1pt} {\kern 1pt} m(\xi ),{\kern 1pt} {\kern 1pt} {\kern 1pt} {%
\kern 1pt} {\kern 1pt} {\kern 1pt} {\kern 1pt} {\kern 1pt} {\kern 1pt} {%
\kern 1pt} {\kern 1pt} {\kern 1pt} f\in BMOA,{\kern 1pt} {\kern 1pt} {\kern %
1pt} {\kern 1pt} {\kern 1pt} {\kern 1pt} {\kern 1pt} {\kern 1pt} {\kern 1pt}
{\kern 1pt} {\kern 1pt} {\kern 1pt} {\kern 1pt} {\kern 1pt} h\in H^{1} .
\end{equation*}

If $h\in H^{1} $ then

\begin{equation*}
\left|\int _{T}\overline{f(r\xi )} {\kern 1pt} {\kern 1pt} {\kern 1pt} {%
\kern 1pt} h(\xi ){\kern 1pt} {\kern 1pt} {\kern 1pt} d{\kern 1pt} {\kern 1pt%
} m(\xi ){\kern 1pt} \right|{\kern 1pt} {\kern 1pt} {\kern 1pt} {\kern 1pt} {%
\kern 1pt} {\kern 1pt} {\kern 1pt} {\kern 1pt} \le {\kern 1pt} {\kern 1pt} {%
\kern 1pt} \sum \nolimits _{n=1}^{\infty }\left|\hat{f}(n){\kern 1pt} \right|%
{\kern 1pt} {\kern 1pt} r^{n} \left|\int _{T}\xi {\kern 1pt} {\kern 1pt}
^{-n} {\kern 1pt} {\kern 1pt} {\kern 1pt} {\kern 1pt} h(\xi ){\kern 1pt} {%
\kern 1pt} {\kern 1pt} d{\kern 1pt} {\kern 1pt} m(\xi )\right|{\kern 1pt} {%
\kern 1pt} {\kern 1pt} {\kern 1pt} ={\kern 1pt}
\end{equation*}

\begin{equation*}
={\kern 1pt} {\kern 1pt} \sum \nolimits _{n=1}^{\infty }\left|\hat{f}(n){%
\kern 1pt} \right|{\kern 1pt} {\kern 1pt} r^{n} \left|\hat{h}(n)\right|{%
\kern 1pt} {\kern 1pt} {\kern 1pt} {\kern 1pt} {\kern 1pt} {\kern 1pt} {%
\kern 1pt} {\kern 1pt} {\kern 1pt} {\kern 1pt} {\kern 1pt} {\kern 1pt} {%
\kern 1pt} \le \left\| \mu \right\| {\kern 1pt} {\kern 1pt} {\kern 1pt} {%
\kern 1pt} {\kern 1pt} \sum \nolimits _{n=1}^{\infty }{\kern 1pt} \frac{1}{n}
{\kern 1pt} \left|\hat{h}(n)\right|{\kern 1pt} r^{n} .
\end{equation*}

By Hardy's inequality

\begin{equation*}
{\kern 1pt} {\kern 1pt} {\kern 1pt} {\kern 1pt} {\kern 1pt} \sum \nolimits
_{n=1}^{\infty }{\kern 1pt} \frac{1}{n} {\kern 1pt} \left|\hat{h}(n)\right|{%
\kern 1pt} {\kern 1pt} {\kern 1pt} {\kern 1pt} {\kern 1pt} \le {\kern 1pt} {%
\kern 1pt} {\kern 1pt} {\kern 1pt} \pi \left\| h\right\| _{H^{1} } .
\end{equation*}

Therefore

\begin{equation*}
\left\| f\right\| _{BMOA} =\sup \left\{\mathop{\lim }\limits_{r\to 1}
\left|\int _{T}\overline{f(r\xi )} {\kern 1pt} {\kern 1pt} {\kern 1pt} {%
\kern 1pt} h(\xi ){\kern 1pt} {\kern 1pt} {\kern 1pt} d{\kern 1pt} {\kern 1pt%
} m(\xi )\right|,{\kern 1pt} {\kern 1pt} {\kern 1pt} {\kern 1pt} {\kern 1pt}
{\kern 1pt} \left\| h\right\| _{H^{1} } \le {\kern 1pt} {\kern 1pt} 1{\kern %
1pt} {\kern 1pt} \right\}{\kern 1pt} {\kern 1pt} {\kern 1pt} \le \pi {\kern %
1pt} {\kern 1pt} {\kern 1pt} \left\| \mu \right\| {\kern 1pt} {\kern 1pt} {%
\kern 1pt} {\kern 1pt} {\kern 1pt} {\kern 1pt} <\infty .{\kern 1pt} {\kern %
1pt} {\kern 1pt}
\end{equation*}

\textbf{Corollary.} $J\subset BMOA$ .

\

\

\noindent {\small Department of Mathematics\newline
Technical University\newline
25, Tsanko Dijstabanov,\newline
Plovdiv, Bulgaria\newline
e-mail: peyyyo@mail.bg}


\begin{thebibliography}{9}
\bibitem{1} J. B. Garnett. \emph{Bounded analytic functions.} Academic.
Academic Press, New York, 1981.

\bibitem{2} T. H. MacGregor. \emph{Analytic and univalent functions with
integral representations involving complex measures.} Indiana Univ. Math.
J., 36, 1987, 109-130.

\bibitem{3} R. A. Hibschweiler, T. H. MacGregor. \emph{Multipliers of
families of Cauchy-Stieltjes transforms.} Trans. Amer. Math. Soc, 331
(1992), 377-394.

\bibitem{4} P. Stoilov. \emph{Some inequalities for analytic functions and
applications to the multipliers of Cauchy-Stieltjcs integrals.}  Mathematics
and Education in Mathematics, Sofia, 1994, 188-193.
\end{thebibliography}
\end{document}